\documentclass[11pt,a4paper]{article}
\usepackage[margin=30mm,marginparwidth=20mm]{geometry}
\usepackage[authoryear,round]{natbib}
\newenvironment{keyword}{\par\small\noindent\textbf{Keywords: }\ignorespaces}{\par\medskip}
\usepackage[T1]{fontenc}
\usepackage[utf8]{inputenc}
\usepackage[british]{babel}
\usepackage{amssymb}
\usepackage{microtype}
\usepackage{xcolor}
\usepackage{graphicx}
\usepackage{subcaption}
\let\cite\citep

\usepackage{array,booktabs}
\usepackage{adjustbox}
\usepackage{longtable}
\usepackage{blkarray}
\usepackage{float}
\usepackage{mathtools}
\usepackage{multirow}
\usepackage{multicol}
\usepackage{pgfplots}
\pgfplotsset{compat=1.18}
\usepackage{pgfplotstable}
\usepackage{pifont}
\usepackage{etoolbox}
\usepackage{siunitx}
\providecommand{\qty}{\SI}%
\providecommand{\unit}{\si}
\DeclareSIUnit{\byte}{B}
\usepackage{tabularx}
\usepackage{tikz}

\robustify\bfseries
\usetikzlibrary{shapes.geometric, positioning, calc}

\newcommand{\N}{\ensuremath{\mathbb{N}}}

\newcommand{\R}{\ensuremath{\mathbb{R}}}

\newcommand{\aC}{\ensuremath{\mathcal{S}}} %

\newcommand{\cfmaxn}{\ensuremath{\bar{k}^{\text{cut}}}}
\newcommand{\cfmaxrel}{\ensuremath{\bar{\alpha}^{\text{cut}}}}
\newcommand{\cfmaxreliter}{\ensuremath{\bar{\alpha}^{\text{cut}}_{\text{iter}}}}
\newcommand{\cfpb}{\ensuremath{\beta^{\text{cut}}}}

\definecolor{rwth-blue}{RGB}{0,84,159}
\definecolor{rwth-lblue}{RGB}{142,186,229}
\definecolor{rwth-purple}{RGB}{122,111,172}
\definecolor{rwth-violet}{RGB}{97,33,88}
\definecolor{rwth-carmine}{RGB}{161,16,53}
\definecolor{rwth-red}{RGB}{204,7,30}
\definecolor{rwth-magenta}{RGB}{227,0,102}
\definecolor{rwth-orange}{RGB}{246,168,0}
\definecolor{rwth-yellow}{RGB}{255,237,0}
\definecolor{rwth-grass}{RGB}{189,171,39}
\definecolor{rwth-green}{RGB}{87,171,39}
\definecolor{rwth-cyan}{RGB}{0,152,161}
\definecolor{rwth-teal}{RGB}{0,97,101}
\definecolor{rwth150red}{HTML}{a11035}

\DeclareDocumentCommand\colorsOn{}{
  \definecolor{visual-red}{RGB}{255,0,0}
  \definecolor{visual-blue}{RGB}{0,0,255}
  \definecolor{visual-green}{RGB}{0,128,0}
  \definecolor{visual-purple}{RGB}{220,0,220}
  \definecolor{visual-orange}{RGB}{255,128,0}
  \definecolor{visual-cyan}{RGB}{0,200,200}
}

\colorsOn

\usepackage[ruled,vlined]{algorithm2e}
\SetKwInput{KwData}{Input}
\SetKwInput{KwResult}{Output}
\SetKw{KwGoto}{goto}
\SetKw{KwWith}{with}
\SetKw{AND}{and}
\SetKw{OR}{or}
\SetKw{GOTO}{go to}
\SetKw{STOP}{STOP.}

\definecolor{ibm-blue}{RGB}{100,143,255}
\definecolor{ibm-purple}{RGB}{120,94,240}
\definecolor{ibm-red}{RGB}{220,38,127}
\definecolor{ibm-orange}{RGB}{254,97,0}
\definecolor{ibm-yellow}{RGB}{255,176,0}

\colorlet{color1}{ibm-blue}
\colorlet{color2}{ibm-orange}
\colorlet{color3}{ibm-red}
\colorlet{color4}{ibm-purple}
\colorlet{color5}{ibm-yellow}

\newlength{\normalparindent}
\AtBeginDocument{\setlength{\normalparindent}{\parindent}}

\usepackage{hyperref}
\usepackage{orcidlink}
\hypersetup{colorlinks=true,allcolors=blue!40!black}
\hypersetup{pdfauthor={Tim Donkiewicz and Oliver Gaul},pdftitle={Benders Cut Filtering}}
\usepackage[capitalise,nameinlink,noabbrev]{cleveref}
\newcommand{\PaperTitle}{Benders Cut Filtering for Affine Potential-Based Flow Problems with Robustness Scenarios and Topology Switching}
\newcommand{\PaperAbstract}{Many large-scale optimization problems decompose into a master problem and scenario subproblems.
Benders decomposition exploits this structure by iteratively solving the subproblems and adding constraints, called Benders cuts, to the master problem.
We study \emph{Benders cut filtering}: given a set of Benders cuts produced in one iteration, which subset should be added to the master problem?
Adding all generated cuts can cause rapid master problem growth, whereas retaining fewer cuts may require more Benders iterations.
To our knowledge, few prior works have studied this question for Benders cuts.
Classical cut filtering techniques require modification because Benders cuts enforce constraints omitted from the master problem, whereas conventional cutting planes are optional inequalities that strengthen the relaxation.
We propose violation-based filtering, diversity-based filtering via $k$-medoids clustering on cosine distances, and a hybrid that selects a most-violated cut per cluster.
We compare these strategies with no filtering, random filtering, and our adaptation of SCIP's hybrid cut selector on an affine potential-based flow problem with topology switching and robustness scenarios.
On instances solved by all configurations, our hybrid approach reduces the shifted geometric mean solve time by about \qty{45}{\percent} compared to no filtering, and by \qty{39}{\percent} compared to our adaptation of SCIP's hybrid cut selector.
While violation-based filtering achieves a similar runtime, the hybrid approach requires about \qty{17}{\percent} fewer callback rounds.
These results demonstrate the effectiveness of Benders cut filtering on this problem class and motivate evaluation in other domains.}
\newcommand{\PaperKeywords}{integer programming; Benders decomposition; cut filtering; potential-based flows; network optimization}

\usepackage{authblk}
\usetikzlibrary{svg.path}
\definecolor{orcidlogocol}{HTML}{A6CE39}
\tikzset{
    orcidlogo/.pic={
        \fill[orcidlogocol] svg{M256,128c0,70.7-57.3,128-128,128C57.3,256,0,198.7,0,128C0,57.3,57.3,0,128,0C198.7,0,256,57.3,256,128z};
        \fill[white] svg{M86.3,186.2H70.9V79.1h15.4v48.4V186.2z}
        svg{M108.9,79.1h41.6c39.6,0,57,28.3,57,53.6c0,27.5-21.5,53.6-56.8,53.6h-41.8V79.1z M124.3,172.4h24.5c34.9,0,42.9-26.5,42.9-39.7c0-21.5-13.7-39.7-43.7-39.7h-23.7V172.4z}
        svg{M88.7,56.8c0,5.5-4.5,10.1-10.1,10.1c-5.6,0-10.1-4.6-10.1-10.1c0-5.6,4.5-10.1,10.1-10.1C84.2,46.7,88.7,51.3,88.7,56.8z};
    }
}

\newcommand{\paperorcid}[1]{\href{https://orcid.org/#1}{%
  \begin{tikzpicture}[xscale=0.025,yscale=-0.025,transform shape]
    \pic{orcidlogo};
  \end{tikzpicture}}}

\begin{document}
\title{\PaperTitle}
\author{Tim Donkiewicz\,\protect\paperorcid{0000-0002-5721-3563}}
\author{Oliver Gaul\,\protect\paperorcid{0000-0002-2131-1911}}
\affil{Chair of Operations Research, RWTH Aachen University\\
  Kackertstra\ss e 7, 52072 Aachen, Germany}
\date{}
\maketitle
\begin{abstract}
\PaperAbstract
\end{abstract}
\begingroup\small
\noindent\textbf{Corresponding author.} Tim Donkiewicz: \href{mailto:donkiewicz@or.rwth-aachen.de}{donkiewicz@or.rwth-aachen.de}.
\par\endgroup\medskip
\section{Introduction}
\label{sec:introduction}
Benders decomposition~\cite{benders1962partitioning} splits an optimization problem into a master problem and one or more subproblems.
The method repeatedly solves the master problem and adds constraints, called \emph{Benders cuts}, obtained by solving the subproblems.
Branch-and-Benders-cut~\cite[cf.][]{rahmaniani2017benders} embeds this principle into branch-and-bound: subproblems are solved in a callback at each feasible master solution, and the resulting Benders cuts are added as lazy constraints.
When the number of subproblems is large, each iteration during branch-and-Benders-cut can produce many such constraints, and adding all of them leads to rapid master problem growth and potentially longer re-optimization times.
Still, the standard approach in the literature is to add all generated Benders cuts to the master problem, either individually or in aggregated form~\cite{rahmaniani2017benders}.

We study a problem that we call \emph{Benders cut filtering}: given the pool of Benders cuts generated in an iteration, which subset should be added as constraints to the master problem?
This is particularly relevant when there are many subproblems that are inexpensive to solve.
Prior work has filtered Benders cuts with machine-learning classifiers trained offline on labeled examples~\cite{hasan2023accelerating, jia2021benders}.
We study filtering rules that do not require offline training on solved instances to determine which Benders cuts to add.

Outside the Benders context, cutting plane selection is well studied in mixed-integer programming (MIP) branch-and-cut solvers~\cite{achterberg2007constraint,andreello2007embedding,turner2023adaptive}.
Solvers manage cut pools and apply quality criteria such as violation, pairwise distance, or objective parallelism to decide which cuts to enforce.
These cutting planes aim to strengthen the relaxation, whereas Benders cuts are \emph{constraints} initially omitted from the master problem.
SCIP adds Benders cuts as constraints by default, so they do not pass through its cut selector~\cite{maher2026scipbenders}.

A separate line of work, \emph{Benders cut selection}, addresses a different problem not considered in this paper: choosing a strongest cut among multiple dual solutions of a single subproblem~\cite{magnanti1981accelerating,kaltis2025literature}.

We propose three filtering strategies for a budget of $k$ cuts:
\emph{Violation-based} filtering retains the $k$ most-violated cuts.
\emph{Diversity-based} filtering applies $k$-medoids clustering on pairwise cosine distances between cut coefficient vectors and selects one cut closest to the centroid of each cluster.
\emph{Hybrid} filtering clusters identically but selects a most-violated cut per cluster, combining violation with diversity.
As baselines, we adapt SCIP's default cut selector and consider random filtering.
We specify the cut budget by a fixed count, a fraction of the number of scenarios, or a fraction of the number of cuts generated in an iteration.

We demonstrate our Benders cut filtering on an affine potential-based flow problem~\cite{pfetsch_potential-based_2026} with topology switching and robustness scenarios.
The master problem contains switching and generation decisions, while each robustness scenario is a single-edge outage in the underlying network graph.
Every such outage is represented by a flow subproblem, solvable using systems of linear equations~\cite{gaul2026efficient}.
Since the subproblems are structurally similar, differing by at most two arcs in a graph typically containing hundreds of arcs, the generated cuts may also be similar, making filtering particularly relevant.

On the jointly solved instances, Hybrid has the lowest shifted geometric mean solve time, closely followed by Violation.
The study shows that Benders cut filtering can substantially improve solve times on suitable problems.
Assessing the generalizability of the approach to other problem classes is left for future work.

\section{Related Work}
\label{sec:related-work}

\citet{rahmaniani2017benders} provide a comprehensive survey of Benders decomposition.
\citet{van1969shaped} combine dual information from multiple subproblems to generate a single aggregated cut.
\citet{birge1988multicut} introduce a multi-cut variant, which adds separate cuts for the individual subproblems.
This can reduce the number of iterations at the cost of a larger master problem.

Cut aggregation requires choosing which subproblems' cuts to combine, a choice that can be adapted during solving~\cite{trukhanov2010adaptive,ramirez2023benders}.
\citet{trukhanov2010adaptive} also test grouping by cut coefficient similarity but find limited runtime benefits due to the cost of comparing coefficients.
\citet{song2015adaptive} adaptively aggregate subproblems instead of their cuts.
\citet{wolf2013dynamic} replace previously added cuts with an aggregated cut, reducing the retained number of cuts.

Other ways to contain master problem growth are partial decomposition and subproblem selection.
Partial decomposition means retaining some subproblem variables in the master problem, representing more of the original problem directly and potentially reducing the number of Benders cuts needed~\cite{crainic2021partial,pauphilet2025random}.
Subproblem selection limits cut generation in each iteration by solving only a subset of subproblems~\cite{donkiewicz2026subproblemselection,blanchot2023benders,bertsimas2025stochastic}.

Existing Benders cut filtering methods include approaches that train classifiers on solved instances to decide which generated cuts to add, using labels based on criteria such as dual bound improvement~\cite{jia2021benders,hasan2023accelerating}.

Cut filtering is well established outside of the Benders context in branch-and-cut solvers, where it is well-known that adding all generated cuts can increase solve times~\cite{andreello2007embedding}.
SCIP's hybrid cut selector combines scoring criteria such as normalized violation with a filter that rejects cuts nearly parallel to those already selected~\cite{gottwald2026sciphybrid}.
Machine learning approaches have also arisen to learn weights for cut scoring that depend on the instance~\cite{turner2023adaptive} or use reinforcement learning to select cuts~\cite{tang2020reinforcement,deza2023machine}.

\paragraph*{Positioning}

Inspired by cut management in MIP solvers, we develop violation-based and clustering-based filtering strategies without offline training, adapted to Benders decomposition.

Our filtering strategies use information from the generated cuts to decide which constraints to add to the master, whereas partial decomposition and subproblem selection act before cut generation.

Our filtering strategies retain selected cuts in their original form, whereas cut aggregation combines information from multiple subproblems into new inequalities.

\section{Problem Description}
\label{sec:potential-based-flows}

We solve problems within the potential-based flow framework described by \citet{pfetsch_potential-based_2026}.
In these network flow problems, potential equations relate arc flows to the potential differences between their endpoints.
Certain arcs can be switched on or off.
An arc that is switched off carries no flow and no longer constrains the potential difference between its endpoints.
The framework covers diverse network optimization problems, including electrical power systems (direct current power flow), gas networks, and traffic assignment.

We solve an extended version with adjustable generation, scenario-based robustness, and overloads as described by~\citet{donkiewicz2026integrated}.
Each scenario corresponds to the failure of a set of arcs.
Switching and generation decisions are made in a first stage before observing which scenario occurs, while flows and potentials adapt per scenario in a second stage.
We minimize switching and generation adjustment costs together with penalties for exceeding arc capacities across all scenarios.

This two-stage structure admits a Benders decomposition~\cite{benders1962partitioning}: the master problem contains switching and generation variables, while each scenario yields a subproblem that checks feasibility and computes overload costs under the given topology and outage.
Bounds on potential differences can make a subproblem infeasible and thus require feasibility cuts.

For the affine potential functions considered here, each subproblem can be solved using a fixed number of systems of linear equations, yielding both primal and dual solutions without solving a linear program~\cite{gaul2026efficient}.

The instances studied in this paper have many independent subproblems, each of which can yield a Benders cut at the same master solution.
Scenarios that share affected regions may produce cuts with similar coefficient vectors, so adding all cuts may grow the master with little bound improvement.
Since the subproblems are inexpensive to solve, the cost of additional subproblem solves may be offset by savings in master solve time, even if filtering requires more Benders iterations.

\section{Benders Cut Filtering Strategies}
\label{sec:cut-filtering}

In each branch-and-Benders-cut callback, we first check connectivity at the current master solution.
If this check yields cuts, we add them as lazy constraints without filtering and return from the callback.
Otherwise, we solve all scenario flow subproblems to obtain a pool $\mathcal{P}$ of violated optimality and feasibility cuts.
We select a subset $\mathcal{C} \subseteq \mathcal{P}$ and add these cuts as lazy constraints.
Correctness and convergence are preserved as long as at least one violated cut is added to reject each master solution with infeasible subproblems or underestimated subproblem costs.

We present three filtering strategies: violation-based scoring, diversity-based clustering, and a hybrid of the two.
Violation-based scoring focuses on the effectiveness of each individual cut, whereas diversity-based clustering aims to avoid adding nearly parallel cuts.
The violation-based and hybrid strategies prioritize feasibility cuts, whereas random and diversity-based filtering do not.

We specify the number $k = |\mathcal{C}| \leq |\mathcal{P}|$ of retained cuts by a fixed count $\cfmaxn \in \N$, a count $\lceil \cfmaxrel \cdot |\aC| \rceil$ based on a fraction $\cfmaxrel \in (0,1]$ of the scenario set $\aC$, or a count $\lceil \cfmaxreliter \cdot |\mathcal{P}| \rceil$ based on a fraction $\cfmaxreliter \in (0,1]$ of the cut pool $\mathcal{P}$.

\subsection{Violation-Based Filtering}
\label{subsec:violation-filtering}

Let $\bar{x}$ denote the current master solution.
For a cut $C_i \in \mathcal{P}$, $i \in [|\mathcal{P}]$, of the form $a_i^\top x \leq b_i$, the \emph{violation} is $v(C_i, \bar{x}) \coloneq \max(0, a_i^\top \bar{x} - b_i)$.
We score optimality cuts by their violation, setting $s(C_i) \coloneq v(C_i, \bar{x})$.
The violation of feasibility cuts depends on their scaling and is, in general, not directly comparable to the violation of optimality cuts.
We thus assign feasibility cuts a score of $s(C_i) \coloneq \infty$, giving them maximum priority over optimality cuts.
Only the $k$ cuts with the largest scores are retained.

For violation-based filtering, we also consider a primal-bound rule with multiplier $\cfpb \geq 1$, as an alternative to specifying a cut count or fraction:
If feasibility cuts exist, we retain all of them and no optimality cuts.
If no feasibility cuts exist and no feasible solution to the original problem is known, we retain all optimality cuts.
With a known feasible solution and no feasibility cuts, we select optimality cuts in decreasing order of violation until their total violation exceeds $\cfpb$ times the nonnegative gap between the best known feasible objective value and the current master objective value or the cut pool is exhausted.

\subsection{Diversity-Based Filtering}
\label{subsec:diversity-kmedoids}

Beyond scoring individual cuts, we also view the entire set of cuts $\mathcal{P}$ together, aiming to retain a diverse subset $\mathcal{C}$.
We compare cuts with nonzero coefficient vectors using the \emph{cosine distance} $d_{\cos}(C_i, C_j) \in [0,2]$ for two cuts $C_i, C_j \in \mathcal{P}$ with coefficient vectors $a_i, a_j \in \R^n$:
\begin{equation}
    d_{\cos}(C_i, C_j) := 1 - \frac{a_i^\top a_j}{\|a_i\|_2 \|a_j\|_2}.
    \label{eq:cosine}
\end{equation}
Cuts whose coefficient vectors point in the same direction have distance $0$.
Cuts whose coefficient vectors point in opposite directions have distance $2$, distinguishing, for example, an upper bound $x_j \leq u$ from a lower bound $x_j \geq \ell$.

We apply $k$-medoids clustering~\cite{kaufman1987clustering} with cosine distance to partition $\mathcal{P}$ into $k$ groups.
For each cluster $\mathcal{K} \subseteq \mathcal{P}$, we then compute the centroid $\bar{a}_{\mathcal{K}} := \frac{1}{|\mathcal{K}|} \sum_{C_j \in \mathcal{K}} a_j$ of its coefficient vectors.
We retain a cut whose coefficient vector is closest to this centroid in Euclidean distance, $C_{\mathcal{K}}^\star \in \arg\min_{C_j \in \mathcal{K}} \| a_j - \bar{a}_{\mathcal{K}} \|_2$.

\subsection{Hybrid Filtering}
\label{subsec:hybrid-filtering}

The hybrid strategy combines violation-based filtering with diversity-based filtering by applying the clustering in \cref{subsec:diversity-kmedoids} and retaining a cut with the largest violation score $s(C)$ from each cluster.

\subsection{SCIP Cut Selection}
\label{subsec:scip-cut-selection}

For comparison, we adapt SCIP's hybrid cut selector~\cite{gottwald2026sciphybrid} to filter Benders constraints.
Here, the score of a cut is a weighted sum of its \emph{efficacy}, its \emph{parallelism} with the objective coefficient vector, and its \emph{integer support} (the fraction of nonzero cut coefficients associated with integer variables).
The efficacy of a cut $C_i$ is $(a_i^\top\bar{x}-b_i)/\max(\|a_i\|_2,\num{1e-6})$.
The parallelism between two cut coefficient vectors is the absolute value of their cosine similarity, $|1-d_{\cos}(C_i,C_j)|$.
Small random perturbations are added to the scores to break ties.
Candidate cuts whose parallelism with a previously selected cut exceeds SCIP's threshold are discarded.

We compute scores from the original Benders cut coefficients, without SCIP's presolve or constraint transformations.
The same selection criteria are applied to feasibility and optimality cuts.
After selection, we discard cuts with efficacy at most $\num{1e-6}$.
We do not replace discarded cuts unless no cut remains, in which case we retain a violated cut of maximum efficacy.

\section{Evaluation}
\label{sec:evaluation}

We evaluate our Benders cut filtering strategies on 131 instances of the affine potential-based flow problem described in \cref{sec:potential-based-flows}.
We first conduct a parameter search and then compare the filtering strategies and baselines.

\paragraph*{Experimental setup}

The instances are generated using the method of \citet{donkiewicz2026instances}, with names referring to the number of nodes.

All experiments run on machines with $2 \times$ Intel Xeon L5630 quad-core processors at \qty{2.13}{\giga\hertz} (8~cores, \qty{16}{\giga\byte} memory).
The implementation is in Julia~1.12.6 using JuMP~\cite{lubin2023jump} with Gurobi~13.0.2~\cite{gurobi} as the MIP solver.
All configurations use four threads for the master problem and a single subproblem thread.

We report arithmetic and shifted geometric mean (SGM) solve times.
For solve times $t_1, \ldots, t_n$ measured in seconds, the SGM is $\bigl(\prod_{i=1}^n (t_i + s)\bigr)^{1/n} - s$, with $s = 10$.

\paragraph*{Configurations}
We compare \textbf{Violation}, \textbf{Diversity}, and \textbf{Hybrid} (\cref{subsec:violation-filtering,subsec:diversity-kmedoids,subsec:hybrid-filtering}) with three baselines.
\textbf{No filtering} adds all cuts, while \textbf{Random} selects $k$ cuts uniformly at random.
For \textbf{SCIP cut selection} (\cref{subsec:scip-cut-selection}), we use the default parameter settings of SCIP~10.0.1, including weights $1$, $0.1$, and $0.1$ for efficacy, objective parallelism, and integer support, respectively.

We also test two variants that allow the number of retained cuts to vary in each iteration: Violation with the primal-bound rule from \cref{subsec:violation-filtering}, and Hybrid with DBSCAN~\cite{ester1996density} in place of $k$-medoids.
Additionally, Hybrid+ supplements the cuts selected by Hybrid with a weighted aggregate of the discarded violated cuts.

\paragraph*{Parameter search}
We conduct the parameter search on 20 instances sampled from the full set.
For Random, Violation, Diversity, Hybrid, and SCIP cut selection, we test fixed cut counts $\cfmaxn \in \{1,2,5,10,20,50\}$ and fractions $\cfmaxrel,\cfmaxreliter \in \{0.01,\allowbreak 0.05,\allowbreak 0.1,\allowbreak 0.25,\allowbreak 0.5,\allowbreak 0.75\}$.
These 90 configurations, SCIP with its default cut limits, and the primal-bound rule with $\cfpb = 1.1$ are compared against no filtering, using a time limit of \qty{1}{\hour}.

The fraction $\cfmaxrel = 0.05$ gives an SGM close to the minimum across Violation, Diversity, and Hybrid and scales the cut budget with the number of scenarios (\cref{fig:gridsearch-heatmap}).
The primal-bound rule does not improve SGM solve time over Violation, Diversity, or Hybrid with $\cfmaxrel = 0.05$ (\cref{fig:gridsearch-gap}).
In preliminary experiments, cosine distance outperforms Jaccard distance.
Replacing $k$-medoids with DBSCAN in Hybrid does not improve runtime.

For the final comparison, we use cosine distance and $k$-medoids for Diversity and Hybrid, $\cfmaxrel = 0.05$ for Random, Violation, Diversity, and Hybrid, and SCIP's default cut limits.
We evaluate all 131 instances with a time limit of \qty{3}{\hour} per instance.

\paragraph*{Results}

\Cref{tab:comparison-time-solve} compares solve times and mean counts on the 126 jointly solved instances, while \cref{fig:performance-profile} includes all 131 instances.
Per-instance results are given in \cref{tab:comparison-time-solve-full}.

\begin{table}[htbp]
    \centering
    \caption{Results on the 126 jointly solved instances, with solved counts for all 131 instances (SGM shift \qty{10}{\second}).}
    \label{tab:comparison-time-solve}
    {\small\begingroup
\setlength{\tabcolsep}{4pt}
\renewcommand{\arraystretch}{1.15}
\sisetup{group-separator={\,},group-minimum-digits=4,detect-weight=true}
\begin{tabular*}{\linewidth}{@{\extracolsep{\fill}}l S[table-format=3.0] *{3}{S[table-format=3.2]} S[table-format=5.2]@{}}
\toprule
 & & \multicolumn{2}{c}{Solve time (\unit{\second})} & \multicolumn{2}{c}{Mean counts} \\
\cmidrule(lr){3-4}\cmidrule(lr){5-6}
Strategy & {Solved} & {Mean} & {SGM} & {\shortstack{Callback\\rounds}} & {\shortstack{Benders\\cuts}} \\
\midrule
No filtering & 126 & 821.27 & 262.75 & \bfseries 302.62 & 34506.23 \\
SCIP cut selection & 131 & 452.05 & 234.70 & 666.86 & \bfseries 1787.89 \\
Random & 131 & 351.05 & 189.78 & 635.26 & 7254.01 \\
\addlinespace[4pt]
Violation & 131 & 249.20 & 144.72 & 435.53 & 4646.98 \\
Diversity & 131 & 242.22 & 153.12 & 395.08 & 4048.74 \\
Hybrid & 131 & \bfseries 225.13 & \bfseries 144.20 & 363.55 & 3672.52 \\
\bottomrule
\end{tabular*}
\endgroup
}
\end{table}

Hybrid has the lowest SGM solve time, closely followed by Violation.
Excluding the 20 instances used in the parameter search leaves the ranking by SGM solve time unchanged and slightly increases Hybrid's advantage over Violation.

Callback rounds in \cref{tab:comparison-time-solve} include connectivity checks, which may reject a candidate before the flow subproblems are solved.
Connectivity cuts are excluded from the cut counts.
Relative to Violation, Hybrid reduces the total number of callback rounds by \qty{16.5}{\percent} and the number of added Benders cuts by \qty{21}{\percent}.
It uses fewer rounds on 97 instances and fewer cuts on 106 instances.
These reductions do not necessarily translate into comparable savings in SGM solve time.

The SCIP selector (reimplemented and adapted to Benders cuts) adds the fewest Benders cuts on average.
However, it requires more callback rounds on average and has higher mean and SGM solve times than Violation and Random.

Hybrid+ does not improve SGM solve time over Hybrid on the jointly solved instances (\cref{tab:comparison-time-solve-full}).

\begin{figure}[htbp]
    \centering
    \input{resources/full_20261023/performance_profile.tikz}
    \caption{Performance profile on all 131 instances. Unsuccessful runs have infinite ratios.}
    \label{fig:performance-profile}
\end{figure}

\section{Summary}
\label{sec:conclusion}

We study Benders cut filtering: selecting which generated cuts to add to the master problem.
We propose and evaluate three strategies that do not require offline training: violation-based filtering, diversity-based filtering using $k$-medoids with cosine distance, and a hybrid that selects a cut with maximum violation score from each cluster.

On the 126 instances solved by all six configurations, Hybrid attains the lowest shifted geometric mean solve time (\qty{144.20}{\second}), closely followed by Violation (\qty{144.72}{\second}), compared with \qty{262.75}{\second} without filtering.
All five filtering strategies solve all 131 instances, whereas the unfiltered method solves 126.
These results support using simple filtering rules when implementing Benders decomposition for this problem class.
In our experiments, strategies using DBSCAN or primal-bound stopping criteria do not improve solve times compared with fixed budgets.
Adding discarded cuts as an aggregate also does not improve over Hybrid.

Further work could investigate other ways of adjusting the number of retained cuts during optimization, by adding fewer cuts while the optimality gap decreases rapidly and more when progress slows.
Cut aging, for example by removing constraints that remain nonbinding over several iterations, could also be explored, although callback interfaces often provide limited control over previously added lazy constraints.
Evaluation in other application domains, such as gas networks or traffic assignment, could help assess whether the observed runtime improvements generalize beyond the problem class studied here.

\paragraph*{Code and data availability}
The implementation, the 131 benchmark instances, and instructions and scripts for reproducing the main comparison are available in
\href{https://github.com/tidonk/TopologyOptimizer.jl/tree/cut_selection-2026-11-10}{TopologyOptimizer.jl, tag \texttt{cut\_selection-2026-11-10}}.

\paragraph*{CRediT authorship contribution statement}
\textbf{Tim Donkiewicz:} Conceptualization, Data curation, Formal analysis, Investigation, Methodology, Software, Validation, Visualization, Writing -- original draft, Writing -- review \& editing.
\textbf{Oliver Gaul:} Software, Writing -- review \& editing.

\paragraph*{Funding}
This research did not receive any specific grant from funding agencies in the public, commercial, or not-for-profit sectors.

\paragraph*{Use of artificial intelligence tools}
The authors used Codex (OpenAI, Astra) and Claude (Anthropic, Opus) through Visual Studio Code to assist with manuscript preparation and editing, the implementation, and the scripts used for analysis and plotting.
All generated code and text were reviewed, verified, and edited by the authors.
The authors take full responsibility for the content of this work.

\bibliographystyle{elsarticle-harv}
\bibliography{literature}

\clearpage

\appendix
\section{Detailed Results}
\label{sec:appendix-results}

\begin{table}[!ht]
    \centering
    \setlength{\abovecaptionskip}{2pt}
    \setlength{\belowcaptionskip}{0pt}
    \caption{Per-instance solve times (\unit{\second}), part 1 of 3. Boldface marks the fastest result. TL: time limit.}
    \label{tab:comparison-time-solve-full}
    \scriptsize
    \setlength{\tabcolsep}{2pt}
    \begin{tabular}{l*{7}{S}}
\toprule
Instance & {No filtering} & {SCIP} & {Random} & {Violation} & {Diversity} & {Hybrid} & {Hybrid+} \\
\midrule
\texttt{85} & 53.08 & 50.90 & 47.95 & \bfseries 44.86 & 51.02 & 50.56 & 49.29 \\
\texttt{92} & \bfseries 36.20 & 42.45 & 39.43 & 37.37 & 42.39 & 43.29 & 42.75 \\
\texttt{92\_2} & 127.10 & 139.43 & 105.25 & \bfseries 70.70 & 84.53 & 77.66 & 76.10 \\
\texttt{92\_3} & 52.13 & 52.89 & 57.20 & \bfseries 45.83 & 52.57 & 48.88 & 49.18 \\
\texttt{94} & 59.91 & 61.71 & 52.49 & \bfseries 52.19 & 56.28 & 54.22 & 54.09 \\
\texttt{98} & \bfseries 45.40 & 69.66 & 60.68 & 57.60 & 56.54 & 54.56 & 54.22 \\
\texttt{100} & 48.25 & 64.60 & 53.20 & \bfseries 43.56 & 47.73 & 50.21 & 50.18 \\
\texttt{105} & \bfseries 43.82 & 65.93 & 51.05 & 45.96 & 53.40 & 51.49 & 51.95 \\
\texttt{110} & 135.20 & 145.16 & 84.98 & \bfseries 80.73 & 85.00 & 82.83 & 84.78 \\
\texttt{112} & 35.87 & 40.87 & 41.09 & \bfseries 35.78 & 46.10 & 41.04 & 41.08 \\
\texttt{114} & \bfseries 38.71 & 43.33 & 39.92 & 41.66 & 47.20 & 46.65 & 47.87 \\
\texttt{115} & 41.19 & 47.88 & 43.26 & \bfseries 38.60 & 46.68 & 46.91 & 46.70 \\
\texttt{116} & 104.74 & 85.02 & 71.64 & \bfseries 56.59 & 70.95 & 67.67 & 69.07 \\
\texttt{117} & \bfseries 36.37 & 55.97 & 39.25 & 37.75 & 41.00 & 42.56 & 42.18 \\
\texttt{120} & 63.16 & 72.68 & 75.47 & 59.76 & 63.50 & 61.51 & \bfseries 59.62 \\
\texttt{120\_2} & 52.63 & 68.92 & 54.25 & \bfseries 48.36 & 56.38 & 57.68 & 53.52 \\
\texttt{122} & 49.17 & 52.87 & 50.05 & \bfseries 45.93 & 52.01 & 49.52 & 50.74 \\
\texttt{122\_2} & 31.84 & 35.47 & \bfseries 30.32 & 31.82 & 37.09 & 36.97 & 37.59 \\
\texttt{124\_2} & 34.21 & 44.33 & \bfseries 33.18 & 34.09 & 39.96 & 39.52 & 39.59 \\
\texttt{125} & 180.72 & 218.89 & 182.48 & 155.67 & 136.15 & \bfseries 127.55 & 136.49 \\
\texttt{126} & 118.97 & 131.94 & 87.56 & \bfseries 80.04 & 85.90 & 87.04 & 91.89 \\
\texttt{128} & 32.47 & 41.37 & 33.35 & \bfseries 31.57 & 38.80 & 37.48 & 37.59 \\
\texttt{128\_2} & \bfseries 40.04 & 71.77 & 49.55 & 51.43 & 56.86 & 51.28 & 52.70 \\
\texttt{130} & 30.55 & 33.62 & 31.48 & 30.38 & 32.51 & \bfseries 29.80 & 33.62 \\
\texttt{131} & 59.95 & 61.52 & 57.29 & \bfseries 50.09 & 60.47 & 56.51 & 56.04 \\
\texttt{132} & 52.26 & 64.73 & 56.17 & \bfseries 47.34 & 55.94 & 51.47 & 55.19 \\
\texttt{134} & 62.93 & 70.79 & 60.34 & 59.58 & 58.62 & 56.49 & \bfseries 56.32 \\
\texttt{138} & 177.04 & 133.33 & 144.22 & \bfseries 100.05 & 120.33 & 114.23 & 120.83 \\
\texttt{138\_1} & \bfseries 31.88 & 46.33 & 37.73 & 35.80 & 39.74 & 39.09 & 43.06 \\
\texttt{146} & 35.84 & 45.32 & 37.12 & \bfseries 34.81 & 40.83 & 41.58 & 42.84 \\
\texttt{147} & 679.67 & 456.44 & 504.30 & 457.82 & 257.91 & \bfseries 211.60 & 222.02 \\
\texttt{148} & 49.94 & 92.31 & 63.68 & \bfseries 49.90 & 54.30 & 52.91 & 52.68 \\
\texttt{148\_1} & 52.36 & 52.70 & 51.68 & \bfseries 48.40 & 52.35 & 53.33 & 58.80 \\
\texttt{152} & \bfseries 64.16 & 104.77 & 101.93 & 65.29 & 74.61 & 68.93 & 73.71 \\
\texttt{153} & 45.16 & 60.28 & 54.76 & \bfseries 44.74 & 52.28 & 48.39 & 50.08 \\
\texttt{155} & 442.78 & 404.81 & 267.21 & 192.34 & 245.22 & \bfseries 186.98 & 211.59 \\
\texttt{155\_2} & 441.84 & 301.23 & 251.90 & \bfseries 164.68 & 219.24 & 195.79 & 201.66 \\
\texttt{156} & 65.12 & 72.93 & 59.98 & \bfseries 51.12 & 60.12 & 58.95 & 57.61 \\
\texttt{159} & 30.87 & 34.44 & 31.03 & 30.94 & \bfseries 30.28 & 30.99 & 33.77 \\
\texttt{159\_2} & 37.89 & \bfseries 36.67 & 38.52 & 38.85 & 38.42 & 37.63 & 39.65 \\
\texttt{161\_1} & 224.73 & 236.14 & 298.03 & 303.86 & 156.88 & 147.44 & \bfseries 143.44 \\
\texttt{162} & 1440.48 & 464.36 & 1069.22 & 304.08 & 391.76 & \bfseries 278.70 & 308.29 \\
\texttt{164} & \bfseries 46.71 & 54.29 & 55.48 & 48.56 & 54.92 & 52.46 & 51.50 \\
\texttt{166\_2} & 57.44 & 89.69 & 58.66 & \bfseries 52.31 & 60.68 & 55.76 & 58.74 \\
\bottomrule
\end{tabular}

\end{table}

\begin{table}[p]
    \centering
    \setlength{\abovecaptionskip}{2pt}
    \setlength{\belowcaptionskip}{0pt}
    \caption{Per-instance solve times (\unit{\second}), part 2 of 3.}
    \label{tab:comparison-time-solve-full-part2}
    \scriptsize
    \setlength{\tabcolsep}{2pt}
    \begin{tabular}{l*{7}{S}}
\toprule
Instance & {No filtering} & {SCIP} & {Random} & {Violation} & {Diversity} & {Hybrid} & {Hybrid+} \\
\midrule
\texttt{167\_2} & 104.23 & 195.92 & 97.74 & 87.41 & \bfseries 73.43 & 91.92 & 77.84 \\
\texttt{169} & 42.84 & 48.57 & 43.74 & \bfseries 42.80 & 46.68 & 46.88 & 48.30 \\
\texttt{169\_1} & 164.55 & 139.20 & 155.97 & \bfseries 76.36 & 111.44 & 99.21 & 100.76 \\
\texttt{174\_1} & 56.63 & 92.05 & 65.84 & \bfseries 55.67 & 66.52 & 62.32 & 58.04 \\
\texttt{179} & \bfseries 39.39 & 60.40 & 48.11 & 45.41 & 47.48 & 47.81 & 49.09 \\
\texttt{180} & 53.79 & 64.86 & 56.31 & \bfseries 49.03 & 52.40 & 54.31 & 51.64 \\
\texttt{187} & 408.90 & 205.28 & 191.18 & \bfseries 137.56 & 158.31 & 144.50 & 147.41 \\
\texttt{188} & 51.48 & 78.20 & 61.04 & \bfseries 47.80 & 53.59 & 55.05 & 55.72 \\
\texttt{191} & 829.13 & 481.98 & 431.86 & 305.56 & 284.60 & 275.95 & \bfseries 271.59 \\
\texttt{191\_1} & 497.80 & 343.15 & 426.51 & \bfseries 189.36 & 248.59 & 190.83 & 224.14 \\
\texttt{191\_2} & 150.08 & 326.81 & 210.54 & 295.01 & 135.12 & 134.21 & \bfseries 131.33 \\
\texttt{194\_3} & 1119.63 & 494.71 & 513.10 & \bfseries 309.62 & 389.35 & 337.31 & 357.96 \\
\texttt{195\_1} & 111.80 & 90.96 & 84.30 & \bfseries 81.07 & 87.88 & 89.55 & 90.66 \\
\texttt{197\_2} & 1751.66 & 774.56 & 930.43 & \bfseries 366.43 & 449.82 & 379.44 & 401.12 \\
\texttt{197\_3} & 1299.01 & 686.76 & 520.30 & \bfseries 334.48 & 361.05 & 354.34 & 374.56 \\
\texttt{207} & \bfseries 51.06 & 95.43 & 71.66 & 69.15 & 69.21 & 60.06 & 63.81 \\
\texttt{211} & 3178.13 & 1137.10 & 814.35 & 699.35 & 716.99 & 693.57 & \bfseries 682.89 \\
\texttt{213} & 1898.62 & 743.94 & 607.86 & \bfseries 274.69 & 495.54 & 414.59 & 427.74 \\
\texttt{215} & 3116.28 & 1550.01 & 1061.61 & \bfseries 573.04 & 666.73 & 600.65 & 617.30 \\
\texttt{217} & 32.51 & 37.21 & 33.01 & \bfseries 31.85 & 37.28 & 38.58 & 39.24 \\
\texttt{218} & 172.12 & 245.11 & 173.39 & 147.45 & \bfseries 136.59 & 148.11 & 145.50 \\
\texttt{218\_1} & 85.05 & 68.75 & 67.27 & \bfseries 57.77 & 59.99 & 64.62 & 66.75 \\
\texttt{220} & 3064.82 & 1343.92 & 1366.01 & 983.84 & \bfseries 782.25 & 845.68 & 912.76 \\
\texttt{222\_2} & 290.81 & 195.72 & 206.89 & \bfseries 128.46 & 172.04 & 146.13 & 153.25 \\
\texttt{223} & 210.48 & 218.78 & 198.48 & \bfseries 112.51 & 146.53 & 134.50 & 134.52 \\
\texttt{232} & 186.74 & 177.81 & 146.38 & \bfseries 105.15 & 114.40 & 109.86 & 112.78 \\
\texttt{233} & 713.89 & 384.54 & 305.86 & 272.18 & 247.50 & 239.62 & \bfseries 238.54 \\
\texttt{234} & 85.60 & 107.65 & 74.88 & \bfseries 66.76 & 81.70 & 78.37 & 75.84 \\
\texttt{235} & 378.04 & 408.69 & 368.70 & 171.92 & 174.36 & \bfseries 166.96 & 178.95 \\
\texttt{237} & 179.32 & 164.66 & \bfseries 94.91 & 97.93 & 113.15 & 106.47 & 99.50 \\
\texttt{239} & 483.82 & 318.50 & 364.65 & \bfseries 201.21 & 252.81 & 214.35 & 228.79 \\
\texttt{240} & 4107.34 & 1229.02 & 1105.57 & 760.58 & 829.52 & 744.93 & \bfseries 718.62 \\
\texttt{242} & 472.81 & 354.64 & 240.08 & 207.42 & 209.11 & \bfseries 159.61 & 231.16 \\
\texttt{248} & 211.55 & 159.71 & 142.98 & 124.01 & 132.92 & \bfseries 122.29 & 136.57 \\
\texttt{250} & {TL} & 2284.75 & 1002.88 & 868.51 & 996.95 & \bfseries 857.80 & 925.27 \\
\texttt{252} & {TL} & 2296.61 & 2000.92 & 1300.75 & 1215.72 & \bfseries 1082.48 & 1206.33 \\
\texttt{252\_2} & 460.24 & 305.00 & 197.91 & \bfseries 189.18 & 192.90 & 204.85 & 197.37 \\
\texttt{255} & 5161.47 & 1520.64 & 1111.52 & \bfseries 778.62 & 867.73 & 844.58 & 838.86 \\
\texttt{258} & 402.99 & 264.10 & 184.38 & \bfseries 157.26 & 171.84 & 183.99 & 169.66 \\
\texttt{258\_1} & 121.20 & 94.15 & 91.51 & \bfseries 76.39 & 84.87 & 77.40 & 77.58 \\
\texttt{262} & 2382.11 & 1146.64 & 924.46 & 642.96 & 679.45 & \bfseries 636.93 & 716.92 \\
\texttt{264} & 221.58 & 236.55 & 164.10 & \bfseries 103.75 & 133.45 & 120.11 & 113.99 \\
\texttt{266} & {TL} & 1346.58 & 1311.57 & 921.59 & 906.45 & \bfseries 804.01 & 869.01 \\
\texttt{269} & 2381.48 & 1366.34 & 1047.32 & 634.13 & 585.11 & 618.68 & \bfseries 564.25 \\
\bottomrule
\end{tabular}

\end{table}

\begin{table}[p]
    \centering
    \setlength{\abovecaptionskip}{2pt}
    \setlength{\belowcaptionskip}{0pt}
    \caption{Per-instance solve times (\unit{\second}), part 3 of 3.}
    \label{tab:comparison-time-solve-full-part3}
    \scriptsize
    \setlength{\tabcolsep}{2pt}
    \begin{tabular}{l*{7}{S}}
\toprule
Instance & {No filtering} & {SCIP} & {Random} & {Violation} & {Diversity} & {Hybrid} & {Hybrid+} \\
\midrule
\texttt{271} & 515.48 & 426.91 & 265.48 & \bfseries 194.00 & 233.37 & 223.95 & 226.94 \\
\texttt{272} & 1077.32 & 726.78 & 439.49 & 310.35 & 368.53 & \bfseries 304.03 & 371.21 \\
\texttt{272\_1} & 1797.61 & 1393.78 & 1205.61 & 809.34 & 627.40 & \bfseries 510.06 & 522.21 \\
\texttt{273} & 2017.78 & 838.43 & 819.10 & 613.54 & 460.07 & 463.56 & \bfseries 452.36 \\
\texttt{277} & 38.52 & 51.38 & 59.98 & \bfseries 38.16 & 42.15 & 52.62 & 49.88 \\
\texttt{284} & 114.16 & 163.10 & \bfseries 91.45 & 95.51 & 92.76 & 92.66 & 93.49 \\
\texttt{293} & 687.82 & 462.00 & 292.17 & \bfseries 260.28 & 274.28 & 263.55 & 281.38 \\
\texttt{294} & 65.10 & 130.79 & 119.73 & 59.52 & 63.00 & \bfseries 56.37 & 66.45 \\
\texttt{303} & 992.79 & 1167.18 & 472.25 & 459.59 & 632.92 & 343.88 & \bfseries 322.92 \\
\texttt{307} & 310.99 & 301.30 & 206.67 & \bfseries 128.62 & 149.62 & 149.82 & 162.78 \\
\texttt{309} & 2756.93 & 1323.38 & 1188.35 & 719.04 & 740.14 & \bfseries 686.20 & 758.05 \\
\texttt{310} & 323.18 & 461.08 & 232.76 & 258.63 & 249.48 & \bfseries 224.41 & 253.69 \\
\texttt{311} & 501.55 & 398.29 & 439.03 & 285.53 & 245.36 & \bfseries 208.56 & 220.08 \\
\texttt{316} & 2936.49 & 1461.01 & 1178.05 & \bfseries 634.85 & 667.60 & 642.63 & 671.37 \\
\texttt{319} & 1459.94 & 915.16 & 669.81 & \bfseries 424.81 & 441.77 & 449.98 & 429.42 \\
\texttt{326} & 364.42 & 223.12 & 190.27 & \bfseries 131.38 & 164.53 & 143.85 & 169.01 \\
\texttt{331\_2} & 2164.01 & 1769.39 & 888.58 & 712.35 & 709.11 & \bfseries 708.76 & 721.57 \\
\texttt{331\_3} & 905.58 & 831.22 & 784.13 & 311.31 & 307.63 & 330.58 & \bfseries 306.33 \\
\texttt{337} & 2187.23 & 1286.05 & 1041.16 & 1623.94 & \bfseries 690.14 & 726.62 & 1127.61 \\
\texttt{341} & 1592.31 & 754.60 & 624.26 & 434.26 & 480.83 & 423.20 & \bfseries 418.25 \\
\texttt{345} & 8175.93 & 1519.26 & 1594.05 & \bfseries 854.47 & 996.45 & 868.31 & 936.17 \\
\texttt{348} & 1107.12 & 525.21 & 478.66 & 341.00 & 399.78 & 338.65 & \bfseries 336.88 \\
\texttt{349} & 2089.09 & 934.88 & 769.15 & \bfseries 526.91 & 593.78 & 572.35 & 611.39 \\
\texttt{349\_1} & 788.31 & 551.11 & 385.22 & \bfseries 262.92 & 348.71 & 297.07 & 314.45 \\
\texttt{352\_2} & 475.25 & 374.37 & 259.89 & \bfseries 192.03 & 226.28 & 209.52 & 212.99 \\
\texttt{353\_2} & 144.29 & 674.33 & 886.76 & \bfseries 129.31 & 299.54 & 148.00 & 284.68 \\
\texttt{356} & 1281.53 & 432.38 & 418.04 & \bfseries 322.04 & 342.24 & 329.12 & 327.51 \\
\texttt{361} & {TL} & 3847.29 & 6473.76 & 3044.78 & 2339.71 & 2147.92 & \bfseries 2115.81 \\
\texttt{361\_2} & \bfseries 122.20 & 1556.70 & 610.47 & 167.25 & 131.29 & 140.36 & 174.17 \\
\texttt{371\_2} & 5819.91 & 3256.20 & 2022.39 & 2065.03 & 1134.93 & \bfseries 1006.98 & 1102.29 \\
\texttt{382} & {TL} & 1780.57 & 1678.57 & 1151.14 & \bfseries 1075.35 & 1079.49 & 1261.70 \\
\texttt{384} & 235.61 & 176.53 & 141.87 & \bfseries 94.98 & 122.49 & 100.30 & 103.25 \\
\texttt{394\_1} & 238.83 & 200.87 & 175.82 & \bfseries 114.62 & 130.13 & 127.61 & 130.76 \\
\texttt{399} & 9133.50 & 2890.10 & 1474.54 & 1002.35 & 861.08 & 808.57 & \bfseries 766.43 \\
\texttt{401} & 446.42 & 387.24 & 253.16 & 202.40 & \bfseries 128.01 & 162.05 & 178.54 \\
\texttt{421\_1} & 747.69 & 822.93 & 420.00 & 334.47 & 342.84 & \bfseries 307.36 & 341.16 \\
\texttt{424\_2} & 1681.06 & 632.39 & 602.45 & 462.54 & \bfseries 441.10 & 473.54 & 459.08 \\
\texttt{431} & 691.96 & 573.67 & 422.05 & \bfseries 242.45 & 308.49 & 297.73 & 325.72 \\
\texttt{465} & 596.95 & \bfseries 396.68 & 431.30 & 533.82 & 442.41 & 470.44 & 455.10 \\
\texttt{494} & 168.87 & 676.94 & 233.06 & \bfseries 141.22 & 153.61 & 155.76 & 157.27 \\
\texttt{518} & 1350.26 & 835.96 & 751.07 & 566.67 & 496.56 & 491.71 & \bfseries 475.11 \\
\texttt{555} & 3008.99 & 1503.63 & 914.28 & 737.03 & \bfseries 678.65 & 690.51 & 730.98 \\
\texttt{605} & \bfseries 45.59 & 372.66 & 62.39 & 71.32 & 74.97 & 49.74 & 50.32 \\
\bottomrule
\end{tabular}

\end{table}

\clearpage
\section{Parameter Search}
The selected instances are \texttt{94}, \texttt{116}, \texttt{117}, \texttt{120}, \texttt{138\_1}, \texttt{146}, \texttt{148\_1}, \texttt{153}, \texttt{162}, \texttt{195\_1}, \texttt{213}, \texttt{215}, \texttt{217}, \texttt{218\_1}, \texttt{223}, \texttt{242}, \texttt{248}, \texttt{353\_2}, \texttt{356}, and \texttt{394\_1}%
.

\begin{figure}[!ht]
    \centering
    \begin{tikzpicture}[x=1cm,y=1cm,font=\scriptsize]
\node[anchor=east] at (-2.15,-.15) {Relative SGM};
\definecolor{gridtile}{HTML}{2166AC}
\filldraw[fill=gridtile,draw=black!25] (-2.0,-.25) rectangle (-1.2,-.05);
\node at (-1.6,-.45) {$\leq 0.5$};
\definecolor{gridtile}{HTML}{A3BFDD}
\filldraw[fill=gridtile,draw=black!25] (-1.2,-.25) rectangle (-0.3999999999999999,-.05);
\node at (-0.7999999999999999,-.45) {0.75};
\definecolor{gridtile}{HTML}{FFFFFF}
\filldraw[fill=gridtile,draw=black!25] (-0.3999999999999999,-.25) rectangle (0.40000000000000013,-.05);
\node at (1.1102230246251565e-16,-.45) {1};
\definecolor{gridtile}{HTML}{D39D6D}
\filldraw[fill=gridtile,draw=black!25] (0.40000000000000036,-.25) rectangle (1.2000000000000004,-.05);
\node at (0.8000000000000004,-.45) {1.5};
\definecolor{gridtile}{HTML}{B35806}
\filldraw[fill=gridtile,draw=black!25] (1.2000000000000002,-.25) rectangle (2.0,-.05);
\node at (1.6,-.45) {$\geq 2$};
\end{tikzpicture}
    \par\medskip
    \begin{subfigure}{\linewidth}
        \centering
        \begin{tikzpicture}[x=1.2cm,y=.52cm,font=\small]
\path[use as bounding box] (-1.55,-.04) rectangle (7.4,5.65);
\node[anchor=east] at (-.12,5.35) {$\cfmaxn$};
\node[font=\scriptsize] at (0.5,5.35) {$1$};
\node[font=\scriptsize] at (1.5,5.35) {$2$};
\node[font=\scriptsize] at (2.5,5.35) {$5$};
\node[font=\scriptsize] at (3.5,5.35) {$10$};
\node[font=\scriptsize] at (4.5,5.35) {$20$};
\node[font=\scriptsize] at (5.5,5.35) {$50$};
\node[anchor=east] at (-.12,4.5) {SCIP};
\definecolor{gridtile}{HTML}{D4A071}
\filldraw[fill=gridtile,draw=white] (0,4) rectangle (1,5);
\node[text=black,font=\footnotesize] at (0.5,4.67) {1.49};
\node[text=black,font=\tiny] at (0.5,4.24) {$S=20$};
\definecolor{gridtile}{HTML}{F7EDE4}
\filldraw[fill=gridtile,draw=white] (1,4) rectangle (2,5);
\node[text=black,font=\footnotesize] at (1.5,4.67) {1.08};
\node[text=black,font=\tiny] at (1.5,4.24) {$S=20$};
\definecolor{gridtile}{HTML}{E1EAF4}
\filldraw[fill=gridtile,draw=white] (2,4) rectangle (3,5);
\node[text=black,font=\footnotesize] at (2.5,4.67) {0.91};
\node[text=black,font=\tiny] at (2.5,4.24) {$S=20$};
\definecolor{gridtile}{HTML}{D5E2EF}
\filldraw[fill=gridtile,draw=white] (3,4) rectangle (4,5);
\node[text=black,font=\footnotesize] at (3.5,4.67) {0.88};
\node[text=black,font=\tiny] at (3.5,4.24) {$S=20$};
\definecolor{gridtile}{HTML}{D5E2EF}
\filldraw[fill=gridtile,draw=white] (4,4) rectangle (5,5);
\node[text=black,font=\footnotesize] at (4.5,4.67) {0.88};
\node[text=black,font=\tiny] at (4.5,4.24) {$S=20$};
\definecolor{gridtile}{HTML}{D3E1EF}
\filldraw[fill=gridtile,draw=white] (5,4) rectangle (6,5);
\node[text=black,font=\footnotesize] at (5.5,4.67) {0.87};
\node[text=black,font=\tiny] at (5.5,4.24) {$S=20$};
\node[anchor=east] at (-.12,3.5) {Random};
\definecolor{gridtile}{HTML}{B35806}
\filldraw[fill=gridtile,draw=white] (0,3) rectangle (1,4);
\node[text=white,font=\footnotesize] at (0.5,3.67) {8.17};
\node[text=white,font=\tiny] at (0.5,3.24) {$S=5$};
\definecolor{gridtile}{HTML}{B35806}
\filldraw[fill=gridtile,draw=white] (1,3) rectangle (2,4);
\node[text=white,font=\footnotesize] at (1.5,3.67) {2.34};
\node[text=white,font=\tiny] at (1.5,3.24) {$S=12$};
\definecolor{gridtile}{HTML}{F6EAE0}
\filldraw[fill=gridtile,draw=white] (2,3) rectangle (3,4);
\node[text=black,font=\footnotesize] at (2.5,3.67) {1.09};
\node[text=black,font=\tiny] at (2.5,3.24) {$S=18$};
\definecolor{gridtile}{HTML}{BBD0E6}
\filldraw[fill=gridtile,draw=white] (3,3) rectangle (4,4);
\node[text=black,font=\footnotesize] at (3.5,3.67) {0.81};
\node[text=black,font=\tiny] at (3.5,3.24) {$S=20$};
\definecolor{gridtile}{HTML}{B1C9E2}
\filldraw[fill=gridtile,draw=white] (4,3) rectangle (5,4);
\node[text=black,font=\footnotesize] at (4.5,3.67) {0.78};
\node[text=black,font=\tiny] at (4.5,3.24) {$S=20$};
\definecolor{gridtile}{HTML}{ADC6E0}
\filldraw[fill=gridtile,draw=white] (5,3) rectangle (6,4);
\node[text=black,font=\footnotesize] at (5.5,3.67) {0.77};
\node[text=black,font=\tiny] at (5.5,3.24) {$S=20$};
\node[anchor=east] at (-.12,2.5) {Violation};
\definecolor{gridtile}{HTML}{C27937}
\filldraw[fill=gridtile,draw=white] (0,2) rectangle (1,3);
\node[text=black,font=\footnotesize] at (0.5,2.67) {1.75};
\node[text=black,font=\tiny] at (0.5,2.24) {$S=17$};
\definecolor{gridtile}{HTML}{F8EFE7}
\filldraw[fill=gridtile,draw=white] (1,2) rectangle (2,3);
\node[text=black,font=\footnotesize] at (1.5,2.67) {1.07};
\node[text=black,font=\tiny] at (1.5,2.24) {$S=19$};
\definecolor{gridtile}{HTML}{79A3CD}
\filldraw[fill=gridtile,draw=white] (2,2) rectangle (3,3);
\node[text=black,font=\footnotesize] at (2.5,2.67) {0.66};
\node[text=black,font=\tiny] at (2.5,2.24) {$S=20$};
\definecolor{gridtile}{HTML}{5086BD}
\filldraw[fill=gridtile,draw=white] (3,2) rectangle (4,3);
\node[text=white,font=\footnotesize] at (3.5,2.67) {0.58};
\node[text=white,font=\tiny] at (3.5,2.24) {$S=20$};
\definecolor{gridtile}{HTML}{4881BB}
\filldraw[fill=gridtile,draw=white] (4,2) rectangle (5,3);
\node[text=white,font=\footnotesize] at (4.5,2.67) {0.57};
\node[text=white,font=\tiny] at (4.5,2.24) {$S=20$};
\draw[black,line width=.8pt] (4.025,2.025) rectangle (4.975,2.975);
\definecolor{gridtile}{HTML}{6998C7}
\filldraw[fill=gridtile,draw=white] (5,2) rectangle (6,3);
\node[text=black,font=\footnotesize] at (5.5,2.67) {0.63};
\node[text=black,font=\tiny] at (5.5,2.24) {$S=20$};
\node[anchor=east] at (-.12,1.5) {Diversity};
\definecolor{gridtile}{HTML}{B35806}
\filldraw[fill=gridtile,draw=white] (0,1) rectangle (1,2);
\node[text=white,font=\footnotesize] at (0.5,1.67) {10.15};
\node[text=white,font=\tiny] at (0.5,1.24) {$S=4$};
\definecolor{gridtile}{HTML}{B35806}
\filldraw[fill=gridtile,draw=white] (1,1) rectangle (2,2);
\node[text=white,font=\footnotesize] at (1.5,1.67) {2.45};
\node[text=white,font=\tiny] at (1.5,1.24) {$S=14$};
\definecolor{gridtile}{HTML}{C4D7E9}
\filldraw[fill=gridtile,draw=white] (2,1) rectangle (3,2);
\node[text=black,font=\footnotesize] at (2.5,1.67) {0.83};
\node[text=black,font=\tiny] at (2.5,1.24) {$S=20$};
\definecolor{gridtile}{HTML}{7BA4CE}
\filldraw[fill=gridtile,draw=white] (3,1) rectangle (4,2);
\node[text=black,font=\footnotesize] at (3.5,1.67) {0.66};
\node[text=black,font=\tiny] at (3.5,1.24) {$S=20$};
\definecolor{gridtile}{HTML}{76A1CC}
\filldraw[fill=gridtile,draw=white] (4,1) rectangle (5,2);
\node[text=black,font=\footnotesize] at (4.5,1.67) {0.65};
\node[text=black,font=\tiny] at (4.5,1.24) {$S=20$};
\definecolor{gridtile}{HTML}{81A8D0}
\filldraw[fill=gridtile,draw=white] (5,1) rectangle (6,2);
\node[text=black,font=\footnotesize] at (5.5,1.67) {0.67};
\node[text=black,font=\tiny] at (5.5,1.24) {$S=20$};
\node[anchor=east] at (-.12,0.5) {Hybrid};
\definecolor{gridtile}{HTML}{B35806}
\filldraw[fill=gridtile,draw=white] (0,0) rectangle (1,1);
\node[text=white,font=\footnotesize] at (0.5,0.67) {2.00};
\node[text=white,font=\tiny] at (0.5,0.24) {$S=17$};
\definecolor{gridtile}{HTML}{D5E2EF}
\filldraw[fill=gridtile,draw=white] (1,0) rectangle (2,1);
\node[text=black,font=\footnotesize] at (1.5,0.67) {0.88};
\node[text=black,font=\tiny] at (1.5,0.24) {$S=20$};
\definecolor{gridtile}{HTML}{86ACD2}
\filldraw[fill=gridtile,draw=white] (2,0) rectangle (3,1);
\node[text=black,font=\footnotesize] at (2.5,0.67) {0.69};
\node[text=black,font=\tiny] at (2.5,0.24) {$S=20$};
\definecolor{gridtile}{HTML}{6998C7}
\filldraw[fill=gridtile,draw=white] (3,0) rectangle (4,1);
\node[text=black,font=\footnotesize] at (3.5,0.67) {0.63};
\node[text=black,font=\tiny] at (3.5,0.24) {$S=20$};
\definecolor{gridtile}{HTML}{6D9BC9}
\filldraw[fill=gridtile,draw=white] (4,0) rectangle (5,1);
\node[text=black,font=\footnotesize] at (4.5,0.67) {0.63};
\node[text=black,font=\tiny] at (4.5,0.24) {$S=20$};
\definecolor{gridtile}{HTML}{78A2CD}
\filldraw[fill=gridtile,draw=white] (5,0) rectangle (6,1);
\node[text=black,font=\footnotesize] at (5.5,0.67) {0.66};
\node[text=black,font=\tiny] at (5.5,0.24) {$S=20$};
\filldraw[fill=white,draw=black!35] (6.2,0) rectangle (7.4,5);
\node[align=center,font=\footnotesize] at (6.8,3.4) {No\\filtering};
\node[font=\footnotesize] at (6.8,2.25) {1.00};
\node[font=\scriptsize] at (6.8,1.25) {$S=20$};
\end{tikzpicture}
        \caption{Fixed cut count.}
        \label{fig:gridsearch-absolute}
    \end{subfigure}
    \par\medskip
    \begin{subfigure}{\linewidth}
        \centering
        \begin{tikzpicture}[x=1.2cm,y=.52cm,font=\small]
\path[use as bounding box] (-1.55,-.04) rectangle (7.4,5.65);
\node[anchor=east] at (-.12,5.35) {$\cfmaxrel$};
\node[font=\scriptsize] at (0.5,5.35) {$0.01$};
\node[font=\scriptsize] at (1.5,5.35) {$0.05$};
\node[font=\scriptsize] at (2.5,5.35) {$0.1$};
\node[font=\scriptsize] at (3.5,5.35) {$0.25$};
\node[font=\scriptsize] at (4.5,5.35) {$0.5$};
\node[font=\scriptsize] at (5.5,5.35) {$0.75$};
\node[anchor=east] at (-.12,4.5) {SCIP};
\definecolor{gridtile}{HTML}{FCFDFE}
\filldraw[fill=gridtile,draw=white] (0,4) rectangle (1,5);
\node[text=black,font=\footnotesize] at (0.5,4.67) {0.99};
\node[text=black,font=\tiny] at (0.5,4.24) {$S=20$};
\definecolor{gridtile}{HTML}{D7E3F0}
\filldraw[fill=gridtile,draw=white] (1,4) rectangle (2,5);
\node[text=black,font=\footnotesize] at (1.5,4.67) {0.88};
\node[text=black,font=\tiny] at (1.5,4.24) {$S=20$};
\definecolor{gridtile}{HTML}{D3E1EE}
\filldraw[fill=gridtile,draw=white] (2,4) rectangle (3,5);
\node[text=black,font=\footnotesize] at (2.5,4.67) {0.87};
\node[text=black,font=\tiny] at (2.5,4.24) {$S=20$};
\draw[black,line width=.8pt] (2.025,4.025) rectangle (2.975,4.975);
\definecolor{gridtile}{HTML}{D5E2EF}
\filldraw[fill=gridtile,draw=white] (3,4) rectangle (4,5);
\node[text=black,font=\footnotesize] at (3.5,4.67) {0.88};
\node[text=black,font=\tiny] at (3.5,4.24) {$S=20$};
\definecolor{gridtile}{HTML}{D7E4F0}
\filldraw[fill=gridtile,draw=white] (4,4) rectangle (5,5);
\node[text=black,font=\footnotesize] at (4.5,4.67) {0.88};
\node[text=black,font=\tiny] at (4.5,4.24) {$S=20$};
\definecolor{gridtile}{HTML}{DAE5F1}
\filldraw[fill=gridtile,draw=white] (5,4) rectangle (6,5);
\node[text=black,font=\footnotesize] at (5.5,4.67) {0.89};
\node[text=black,font=\tiny] at (5.5,4.24) {$S=20$};
\node[anchor=east] at (-.12,3.5) {Random};
\definecolor{gridtile}{HTML}{B35806}
\filldraw[fill=gridtile,draw=white] (0,3) rectangle (1,4);
\node[text=white,font=\footnotesize] at (0.5,3.67) {2.39};
\node[text=white,font=\tiny] at (0.5,3.24) {$S=14$};
\definecolor{gridtile}{HTML}{C2D5E8}
\filldraw[fill=gridtile,draw=white] (1,3) rectangle (2,4);
\node[text=black,font=\footnotesize] at (1.5,3.67) {0.83};
\node[text=black,font=\tiny] at (1.5,3.24) {$S=20$};
\definecolor{gridtile}{HTML}{9CBBDA}
\filldraw[fill=gridtile,draw=white] (2,3) rectangle (3,4);
\node[text=black,font=\footnotesize] at (2.5,3.67) {0.73};
\node[text=black,font=\tiny] at (2.5,3.24) {$S=20$};
\draw[black,line width=.8pt] (2.025,3.025) rectangle (2.975,3.975);
\definecolor{gridtile}{HTML}{A5C1DE}
\filldraw[fill=gridtile,draw=white] (3,3) rectangle (4,4);
\node[text=black,font=\footnotesize] at (3.5,3.67) {0.76};
\node[text=black,font=\tiny] at (3.5,3.24) {$S=20$};
\definecolor{gridtile}{HTML}{D4E1EF}
\filldraw[fill=gridtile,draw=white] (4,3) rectangle (5,4);
\node[text=black,font=\footnotesize] at (4.5,3.67) {0.87};
\node[text=black,font=\tiny] at (4.5,3.24) {$S=20$};
\definecolor{gridtile}{HTML}{EAF0F7}
\filldraw[fill=gridtile,draw=white] (5,3) rectangle (6,4);
\node[text=black,font=\footnotesize] at (5.5,3.67) {0.94};
\node[text=black,font=\tiny] at (5.5,3.24) {$S=20$};
\node[anchor=east] at (-.12,2.5) {Violation};
\definecolor{gridtile}{HTML}{EAF0F7}
\filldraw[fill=gridtile,draw=white] (0,2) rectangle (1,3);
\node[text=black,font=\footnotesize] at (0.5,2.67) {0.94};
\node[text=black,font=\tiny] at (0.5,2.24) {$S=19$};
\definecolor{gridtile}{HTML}{4A82BB}
\filldraw[fill=gridtile,draw=white] (1,2) rectangle (2,3);
\node[text=white,font=\footnotesize] at (1.5,2.67) {0.57};
\node[text=white,font=\tiny] at (1.5,2.24) {$S=20$};
\definecolor{gridtile}{HTML}{588CC0}
\filldraw[fill=gridtile,draw=white] (2,2) rectangle (3,3);
\node[text=black,font=\footnotesize] at (2.5,2.67) {0.59};
\node[text=black,font=\tiny] at (2.5,2.24) {$S=20$};
\definecolor{gridtile}{HTML}{7DA5CE}
\filldraw[fill=gridtile,draw=white] (3,2) rectangle (4,3);
\node[text=black,font=\footnotesize] at (3.5,2.67) {0.67};
\node[text=black,font=\tiny] at (3.5,2.24) {$S=20$};
\definecolor{gridtile}{HTML}{B9CFE5}
\filldraw[fill=gridtile,draw=white] (4,2) rectangle (5,3);
\node[text=black,font=\footnotesize] at (4.5,2.67) {0.80};
\node[text=black,font=\tiny] at (4.5,2.24) {$S=20$};
\definecolor{gridtile}{HTML}{E8EFF6}
\filldraw[fill=gridtile,draw=white] (5,2) rectangle (6,3);
\node[text=black,font=\footnotesize] at (5.5,2.67) {0.93};
\node[text=black,font=\tiny] at (5.5,2.24) {$S=20$};
\node[anchor=east] at (-.12,1.5) {Diversity};
\definecolor{gridtile}{HTML}{E4C3A5}
\filldraw[fill=gridtile,draw=white] (0,1) rectangle (1,2);
\node[text=black,font=\footnotesize] at (0.5,1.67) {1.28};
\node[text=black,font=\tiny] at (0.5,1.24) {$S=19$};
\definecolor{gridtile}{HTML}{83AAD1}
\filldraw[fill=gridtile,draw=white] (1,1) rectangle (2,2);
\node[text=black,font=\footnotesize] at (1.5,1.67) {0.68};
\node[text=black,font=\tiny] at (1.5,1.24) {$S=20$};
\definecolor{gridtile}{HTML}{749FCB}
\filldraw[fill=gridtile,draw=white] (2,1) rectangle (3,2);
\node[text=black,font=\footnotesize] at (2.5,1.67) {0.65};
\node[text=black,font=\tiny] at (2.5,1.24) {$S=20$};
\draw[black,line width=.8pt] (2.025,1.025) rectangle (2.975,1.975);
\definecolor{gridtile}{HTML}{89AED3}
\filldraw[fill=gridtile,draw=white] (3,1) rectangle (4,2);
\node[text=black,font=\footnotesize] at (3.5,1.67) {0.69};
\node[text=black,font=\tiny] at (3.5,1.24) {$S=20$};
\definecolor{gridtile}{HTML}{C9DAEB}
\filldraw[fill=gridtile,draw=white] (4,1) rectangle (5,2);
\node[text=black,font=\footnotesize] at (4.5,1.67) {0.84};
\node[text=black,font=\tiny] at (4.5,1.24) {$S=20$};
\definecolor{gridtile}{HTML}{F7FAFC}
\filldraw[fill=gridtile,draw=white] (5,1) rectangle (6,2);
\node[text=black,font=\footnotesize] at (5.5,1.67) {0.98};
\node[text=black,font=\tiny] at (5.5,1.24) {$S=20$};
\node[anchor=east] at (-.12,0.5) {Hybrid};
\definecolor{gridtile}{HTML}{A9C4DF}
\filldraw[fill=gridtile,draw=white] (0,0) rectangle (1,1);
\node[text=black,font=\footnotesize] at (0.5,0.67) {0.76};
\node[text=black,font=\tiny] at (0.5,0.24) {$S=20$};
\definecolor{gridtile}{HTML}{6494C5}
\filldraw[fill=gridtile,draw=white] (1,0) rectangle (2,1);
\node[text=black,font=\footnotesize] at (1.5,0.67) {0.62};
\node[text=black,font=\tiny] at (1.5,0.24) {$S=20$};
\draw[black,line width=.8pt] (1.025,0.025) rectangle (1.975,0.975);
\definecolor{gridtile}{HTML}{6D9AC8}
\filldraw[fill=gridtile,draw=white] (2,0) rectangle (3,1);
\node[text=black,font=\footnotesize] at (2.5,0.67) {0.63};
\node[text=black,font=\tiny] at (2.5,0.24) {$S=20$};
\definecolor{gridtile}{HTML}{8FB2D5}
\filldraw[fill=gridtile,draw=white] (3,0) rectangle (4,1);
\node[text=black,font=\footnotesize] at (3.5,0.67) {0.71};
\node[text=black,font=\tiny] at (3.5,0.24) {$S=20$};
\definecolor{gridtile}{HTML}{CADAEB}
\filldraw[fill=gridtile,draw=white] (4,0) rectangle (5,1);
\node[text=black,font=\footnotesize] at (4.5,0.67) {0.85};
\node[text=black,font=\tiny] at (4.5,0.24) {$S=20$};
\definecolor{gridtile}{HTML}{EAF0F7}
\filldraw[fill=gridtile,draw=white] (5,0) rectangle (6,1);
\node[text=black,font=\footnotesize] at (5.5,0.67) {0.94};
\node[text=black,font=\tiny] at (5.5,0.24) {$S=20$};
\filldraw[fill=white,draw=black!35] (6.2,0) rectangle (7.4,5);
\node[align=center,font=\footnotesize] at (6.8,3.4) {No\\filtering};
\node[font=\footnotesize] at (6.8,2.25) {1.00};
\node[font=\scriptsize] at (6.8,1.25) {$S=20$};
\end{tikzpicture}
        \caption{Fraction of all scenarios.}
        \label{fig:gridsearch-scenarios}
    \end{subfigure}
    \par\medskip
    \begin{subfigure}{\linewidth}
        \centering
        \begin{tikzpicture}[x=1.2cm,y=.52cm,font=\small]
\path[use as bounding box] (-1.55,-.04) rectangle (7.4,5.65);
\node[anchor=east] at (-.12,5.35) {$\cfmaxreliter$};
\node[font=\scriptsize] at (0.5,5.35) {$0.01$};
\node[font=\scriptsize] at (1.5,5.35) {$0.05$};
\node[font=\scriptsize] at (2.5,5.35) {$0.1$};
\node[font=\scriptsize] at (3.5,5.35) {$0.25$};
\node[font=\scriptsize] at (4.5,5.35) {$0.5$};
\node[font=\scriptsize] at (5.5,5.35) {$0.75$};
\node[anchor=east] at (-.12,4.5) {SCIP};
\definecolor{gridtile}{HTML}{EED9C6}
\filldraw[fill=gridtile,draw=white] (0,4) rectangle (1,5);
\node[text=black,font=\footnotesize] at (0.5,4.67) {1.17};
\node[text=black,font=\tiny] at (0.5,4.24) {$S=20$};
\definecolor{gridtile}{HTML}{FCF9F6}
\filldraw[fill=gridtile,draw=white] (1,4) rectangle (2,5);
\node[text=black,font=\footnotesize] at (1.5,4.67) {1.03};
\node[text=black,font=\tiny] at (1.5,4.24) {$S=20$};
\definecolor{gridtile}{HTML}{FDFDFE}
\filldraw[fill=gridtile,draw=white] (2,4) rectangle (3,5);
\node[text=black,font=\footnotesize] at (2.5,4.67) {0.99};
\node[text=black,font=\tiny] at (2.5,4.24) {$S=20$};
\definecolor{gridtile}{HTML}{E3EBF4}
\filldraw[fill=gridtile,draw=white] (3,4) rectangle (4,5);
\node[text=black,font=\footnotesize] at (3.5,4.67) {0.92};
\node[text=black,font=\tiny] at (3.5,4.24) {$S=20$};
\definecolor{gridtile}{HTML}{D6E3F0}
\filldraw[fill=gridtile,draw=white] (4,4) rectangle (5,5);
\node[text=black,font=\footnotesize] at (4.5,4.67) {0.88};
\node[text=black,font=\tiny] at (4.5,4.24) {$S=20$};
\definecolor{gridtile}{HTML}{E1EAF4}
\filldraw[fill=gridtile,draw=white] (5,4) rectangle (6,5);
\node[text=black,font=\footnotesize] at (5.5,4.67) {0.91};
\node[text=black,font=\tiny] at (5.5,4.24) {$S=20$};
\node[anchor=east] at (-.12,3.5) {Random};
\definecolor{gridtile}{HTML}{B35806}
\filldraw[fill=gridtile,draw=white] (0,3) rectangle (1,4);
\node[text=white,font=\footnotesize] at (0.5,3.67) {4.54};
\node[text=white,font=\tiny] at (0.5,3.24) {$S=8$};
\definecolor{gridtile}{HTML}{B35806}
\filldraw[fill=gridtile,draw=white] (1,3) rectangle (2,4);
\node[text=white,font=\footnotesize] at (1.5,3.67) {2.25};
\node[text=white,font=\tiny] at (1.5,3.24) {$S=13$};
\definecolor{gridtile}{HTML}{D8AA80}
\filldraw[fill=gridtile,draw=white] (2,3) rectangle (3,4);
\node[text=black,font=\footnotesize] at (2.5,3.67) {1.43};
\node[text=black,font=\tiny] at (2.5,3.24) {$S=17$};
\definecolor{gridtile}{HTML}{F8F0E8}
\filldraw[fill=gridtile,draw=white] (3,3) rectangle (4,4);
\node[text=black,font=\footnotesize] at (3.5,3.67) {1.07};
\node[text=black,font=\tiny] at (3.5,3.24) {$S=19$};
\definecolor{gridtile}{HTML}{FDFBF9}
\filldraw[fill=gridtile,draw=white] (4,3) rectangle (5,4);
\node[text=black,font=\footnotesize] at (4.5,3.67) {1.02};
\node[text=black,font=\tiny] at (4.5,3.24) {$S=19$};
\definecolor{gridtile}{HTML}{ECF2F8}
\filldraw[fill=gridtile,draw=white] (5,3) rectangle (6,4);
\node[text=black,font=\footnotesize] at (5.5,3.67) {0.94};
\node[text=black,font=\tiny] at (5.5,3.24) {$S=20$};
\node[anchor=east] at (-.12,2.5) {Violation};
\definecolor{gridtile}{HTML}{F7EDE4}
\filldraw[fill=gridtile,draw=white] (0,2) rectangle (1,3);
\node[text=black,font=\footnotesize] at (0.5,2.67) {1.08};
\node[text=black,font=\tiny] at (0.5,2.24) {$S=19$};
\definecolor{gridtile}{HTML}{94B5D7}
\filldraw[fill=gridtile,draw=white] (1,2) rectangle (2,3);
\node[text=black,font=\footnotesize] at (1.5,2.67) {0.71};
\node[text=black,font=\tiny] at (1.5,2.24) {$S=20$};
\definecolor{gridtile}{HTML}{7BA4CE}
\filldraw[fill=gridtile,draw=white] (2,2) rectangle (3,3);
\node[text=black,font=\footnotesize] at (2.5,2.67) {0.66};
\node[text=black,font=\tiny] at (2.5,2.24) {$S=20$};
\definecolor{gridtile}{HTML}{8DB0D4}
\filldraw[fill=gridtile,draw=white] (3,2) rectangle (4,3);
\node[text=black,font=\footnotesize] at (3.5,2.67) {0.70};
\node[text=black,font=\tiny] at (3.5,2.24) {$S=20$};
\definecolor{gridtile}{HTML}{B5CCE3}
\filldraw[fill=gridtile,draw=white] (4,2) rectangle (5,3);
\node[text=black,font=\footnotesize] at (4.5,2.67) {0.79};
\node[text=black,font=\tiny] at (4.5,2.24) {$S=20$};
\definecolor{gridtile}{HTML}{E9F0F7}
\filldraw[fill=gridtile,draw=white] (5,2) rectangle (6,3);
\node[text=black,font=\footnotesize] at (5.5,2.67) {0.93};
\node[text=black,font=\tiny] at (5.5,2.24) {$S=20$};
\node[anchor=east] at (-.12,1.5) {Diversity};
\definecolor{gridtile}{HTML}{B35806}
\filldraw[fill=gridtile,draw=white] (0,1) rectangle (1,2);
\node[text=white,font=\footnotesize] at (0.5,1.67) {12.94};
\node[text=white,font=\tiny] at (0.5,1.24) {$S=3$};
\definecolor{gridtile}{HTML}{B35806}
\filldraw[fill=gridtile,draw=white] (1,1) rectangle (2,2);
\node[text=white,font=\footnotesize] at (1.5,1.67) {2.02};
\node[text=white,font=\tiny] at (1.5,1.24) {$S=14$};
\definecolor{gridtile}{HTML}{C37A39}
\filldraw[fill=gridtile,draw=white] (2,1) rectangle (3,2);
\node[text=black,font=\footnotesize] at (2.5,1.67) {1.73};
\node[text=black,font=\tiny] at (2.5,1.24) {$S=16$};
\definecolor{gridtile}{HTML}{F1F5FA}
\filldraw[fill=gridtile,draw=white] (3,1) rectangle (4,2);
\node[text=black,font=\footnotesize] at (3.5,1.67) {0.96};
\node[text=black,font=\tiny] at (3.5,1.24) {$S=19$};
\definecolor{gridtile}{HTML}{FAF4EE}
\filldraw[fill=gridtile,draw=white] (4,1) rectangle (5,2);
\node[text=black,font=\footnotesize] at (4.5,1.67) {1.05};
\node[text=black,font=\tiny] at (4.5,1.24) {$S=19$};
\definecolor{gridtile}{HTML}{F4F8FB}
\filldraw[fill=gridtile,draw=white] (5,1) rectangle (6,2);
\node[text=black,font=\footnotesize] at (5.5,1.67) {0.97};
\node[text=black,font=\tiny] at (5.5,1.24) {$S=20$};
\node[anchor=east] at (-.12,0.5) {Hybrid};
\definecolor{gridtile}{HTML}{EBF1F8}
\filldraw[fill=gridtile,draw=white] (0,0) rectangle (1,1);
\node[text=black,font=\footnotesize] at (0.5,0.67) {0.94};
\node[text=black,font=\tiny] at (0.5,0.24) {$S=20$};
\definecolor{gridtile}{HTML}{B3CBE3}
\filldraw[fill=gridtile,draw=white] (1,0) rectangle (2,1);
\node[text=black,font=\footnotesize] at (1.5,0.67) {0.79};
\node[text=black,font=\tiny] at (1.5,0.24) {$S=20$};
\definecolor{gridtile}{HTML}{95B6D8}
\filldraw[fill=gridtile,draw=white] (2,0) rectangle (3,1);
\node[text=black,font=\footnotesize] at (2.5,0.67) {0.72};
\node[text=black,font=\tiny] at (2.5,0.24) {$S=20$};
\definecolor{gridtile}{HTML}{A9C4DF}
\filldraw[fill=gridtile,draw=white] (3,0) rectangle (4,1);
\node[text=black,font=\footnotesize] at (3.5,0.67) {0.77};
\node[text=black,font=\tiny] at (3.5,0.24) {$S=20$};
\definecolor{gridtile}{HTML}{C9DAEB}
\filldraw[fill=gridtile,draw=white] (4,0) rectangle (5,1);
\node[text=black,font=\footnotesize] at (4.5,0.67) {0.85};
\node[text=black,font=\tiny] at (4.5,0.24) {$S=20$};
\definecolor{gridtile}{HTML}{E9F0F7}
\filldraw[fill=gridtile,draw=white] (5,0) rectangle (6,1);
\node[text=black,font=\footnotesize] at (5.5,0.67) {0.93};
\node[text=black,font=\tiny] at (5.5,0.24) {$S=20$};
\filldraw[fill=white,draw=black!35] (6.2,0) rectangle (7.4,5);
\node[align=center,font=\footnotesize] at (6.8,3.4) {No\\filtering};
\node[font=\footnotesize] at (6.8,2.25) {1.00};
\node[font=\scriptsize] at (6.8,1.25) {$S=20$};
\end{tikzpicture}
        \caption{Fraction of violated cuts in the current iteration.}
        \label{fig:gridsearch-current}
    \end{subfigure}
    \par\medskip
    \begin{subfigure}{.48\linewidth}
        \centering
        \begin{tikzpicture}[x=1.2cm,y=.52cm,font=\small]
\path[use as bounding box] (-1.55,-.04) rectangle (1,1.65);
\node[anchor=east] at (-.12,1.35) {$\cfpb$};
\node[font=\scriptsize] at (.5,1.35) {$1.1$};
\node[anchor=east] at (-.12,.5) {Violation};
\definecolor{gridtile}{HTML}{BED2E7}
\filldraw[fill=gridtile,draw=white] (0,0) rectangle (1,1);
\node[text=black,font=\footnotesize] at (.5,.67) {0.82};
\node[text=black,font=\tiny] at (.5,.24) {$S=20$};
\end{tikzpicture}
        \caption{Primal bound threshold.}
        \label{fig:gridsearch-gap}
    \end{subfigure}
    \hfill
    \begin{subfigure}{.48\linewidth}
        \centering
        \begin{tikzpicture}[x=1.2cm,y=.52cm,font=\small]
\path[use as bounding box] (-1.55,-.04) rectangle (1,1.65);
\node[font=\scriptsize] at (.5,1.35) {Default};
\node[anchor=east] at (-.12,.5) {SCIP};
\definecolor{gridtile}{HTML}{D9E5F1}
\filldraw[fill=gridtile,draw=white] (0,0) rectangle (1,1);
\node[text=black,font=\footnotesize] at (.5,.67) {0.89};
\node[text=black,font=\tiny] at (.5,.24) {$S=20$};
\end{tikzpicture}
        \caption{SCIP default cut limits.}
        \label{fig:gridsearch-scip-native}
    \end{subfigure}
    \caption{Grid search on 20 instances: shifted geometric mean solve time relative to no filtering. The shift is \qty{10}{\second}, and unsuccessful runs receive \qty{3600}{\second}. $S$ denotes the number of solved instances. Outlines mark the best count or fraction budget for each strategy.}
    \label{fig:gridsearch-heatmap}
\end{figure}

\end{document}